\documentclass[aos]{imsart}

\RequirePackage[OT1]{fontenc}
\RequirePackage{amsthm,amsmath}

\startlocaldefs
\numberwithin{equation}{section}
\theoremstyle{plain}
\newtheorem{thm}{Theorem}

\newtheorem{lem}{Lemma}

\newcommand{\cA}{{\cal A}}

\newcommand{\cN}{{\cal N}}
\newcommand{\cK}{{\cal K}}

\newcommand{\wht}{\widehat}
\newcommand{\wtd}{\widetilde}

\newcommand{\wbox}{\sqcap\llap{$\sqcup$}}

\endlocaldefs

\begin{document}
\begin{frontmatter}



\title{OPTIMAL SEQUENTIAL DETECTION IN MULTI-STREAM DATA} 
\runtitle{OPTIMAL SEQUENTIAL DETECION IN MULTI-STREAM DATA}


\author{\fnms{Hock Peng} \snm{Chan}\ead[label=e1]{stachp@nus.edu.sg}\thanksref{t1}}
\thankstext{t1}{Suported by the National University of Singapore grant R-155-000-158-112}
\address{Department of Statistics\\ and Applied Probability\\ 6 Science Drive 2 \\Singapore 117546}
\affiliation{National University of Singapore}

\runauthor{Hock Peng Chan}



\begin{abstract}
Consider a large number of detectors each generating a data stream. 
The task is to detect online, 
distribution changes in a small fraction of the data streams. 
Previous approaches to this problem include the use of mixture likelihood ratios and sum of CUSUMs.
We provide here extensions and modifications of these approaches that are optimal in detecting normal mean shifts.
We show how the (optimal) detection delay depends on the fraction of data streams undergoing distribution changes as the number of detectors goes to infinity.
There are three detection domains.
In the first domain for moderately large fractions,
immediate detection is possible.
In the second domain for smaller fractions,
the detection delay grows logarithmically with the number of detectors,
with an asymptotic constant extending those in sparse normal mixture detection.
In the third domain for even smaller fractions,
the detection delay lies in the framework of the classical detection delay formula of Lorden.
We show that the optimal detection delay is achieved by the sum of detectability score transformations of either the partial scores or CUSUM scores of the data streams.
\end{abstract}
\end{frontmatter}


\section{Introduction}

Consider $N$ data streams with $X_{nt}$ the observation of the $n$th data stream at time $t$.
We want to detect as quickly as we can a possible change-point $\nu \geq 1$, 
such that for some $\cN \subset \{ 1,\ldots, N \}$,
the post-change observations $X_{nt}$ for $n \in \cN$ (and $t \geq \nu$) have distributions different from the pre-change observations.
Applications for this multi-stream sequential change-point detection problem include hospital management, infectious-disease modeling and target detection.

Tartakovsky and Veervallli \cite{TV08} consider distributed decision-making and optimal fusion,
with minimax,
uniform and Bayesian formulations for sequential detection in multi-stream data.
Though optimal detection is achieved, 
the asymptotics involve $N$ fixed as the average run lengths go to infinity. 

Mei \cite{Mei10} considers distribution changes that do not affect all data streams,
and recommends a sum of CUSUM approach. 
The advantages of his approach are that the distribution changes are not assumed to have occurred simultaneously, 
and the efficient computation of his stopping rule.
However as has been shown in an earlier simulation study, 
the detection delay is relatively large when $\# \cN$, 
the number of data streams undergoing change, 
is small.

Xie and Siegmund \cite{XS13} are the first to look from the perspective of $\# \cN$ small.
They suggest a mixture likelihood ratio (MLR) approach and show via simulation studies the superiority of their MLR stopping rules in detecting over a wide range of $\# \cN$, 
compared to other known approaches.
They also provide analytical approximations to average run lengths and detection delays of their stopping rules that are accurate and useful.
However they do not give any small or moderate $\# \cN$ optimality theory. 

In parallel developments,
motivated by applications in DNA copy-number samples, 
there have been advances made,
see Siegmund, Yakir and Zhang \cite{SYZ11}, Jeng, Cai and Li \cite{JCL13} and Chan and Walther \cite{CW15},
on fixed-sample change-point detection in multiple sequences having a common location index.
The work here also has connections with detection on spatial indices, 
see \cite{ADH05, ADH06, Cha09}. 

In this paper we show that subject to an average run length constraint, 
a modified version of the MLR stopping rule achieves minimum detection delay,
extending the classical single-stream optimal detection of Lorden \cite{Lor71}, 
Pollak \cite{Pol85, Pol87} and Moustakides \cite{Mou86} to multiple data streams,
in the detection of normal mean shifts.
In Section 2 we provide the asymptotic lower bounds of the detection delays for different domains of $\cN$. 
Under the first domain for large $\# \cN$,  
the lower bound is trivially given by 1.
Under the second domain for moderate $\# \cN$,
the lower bound grows logarithmically with $N$.
Under the third domain for small $\# \cN$,
the detection delay grows polynomially with $N$.
In Section 3 we show that a MLR stopping rule that tests against the limits of detectability achieves optimal detection on all three domains.
A window-limited rule, 
suggested in Lai \cite{Lai95}, 
is incorporated into the stopping rule for computational savings.
In Section 4 a numerical study is performed to provide justification for using the MLR stopping rule for finite $N$.
In Section 5 we extend the idea of testing against the limits of detectability on Mei's sum of CUSUM test.
Rather than summing the CUSUM scores as in Mei \cite{Mei10}, 
we suggest instead to sum the detectability score transformations of the CUSUM scores.
Optimality of this procedure is shown but it occurs only when we select the assumed mean shift at a specific value between one to two times the true mean shift,
surprisingly not at the true mean shift itself. 
In Sections 6--8 we provide the proofs of Theorems~\ref{thm1}--\ref{thm3}.

\section{Detection delay lower bound}

Let $X_{nt}$, 
$1 \leq n \leq N$,
$t \geq 1$, 
be distributed as independent N($\mu_{nt},1$). 
Assume that at some unknown time $\nu \geq 1$,
there are mean shifts in a subset $\cN$ of the data streams.
More specifically we assume that 
\begin{equation} \label{munt}
\mu_{nt} = \mu {\bf I}_{\{ t \geq \nu, n \in \cN \}} \mbox{ for some } \mu > 0,
\end{equation}
with ${\bf I}_{\{ 1 \in \cN \}}, \ldots, {\bf I}_{\{ N \in \cN \}}$ i.i.d. Bernoulli($p$) for some $0 < p < 1$.
We shall let $P_\nu$ ($E_\nu$) denote probability measure (expectation) with respect to distribution changes at time $\nu$,
with $\nu=\infty$ indicating no change.
In Appendix~B we provide an analogue of Theorem \ref{thm1} below on a minimax formulation of the problem,
with a constraint on $\sum_{n=1}^N {\bf I}_{\{ n \in \cN \}}$ instead of assuming ${\bf I}_{\{ n \in \cN \}}$ to be i.i.d. Bernoulli.

A standard measure of the performance of a stopping rule $T$,  
see Pollak \cite{Pol85, Pol87},
is the (expected) detection delay 
\begin{equation} \label{DNT}
D_N(T) := \sup_{1 \leq \nu < \infty} E_\nu (T-\nu+1|T \geq \nu),
\end{equation}
subject to the constraint that ARL($T$) (:=$E_\infty T$) $\geq \gamma$ for some $\gamma \geq 1$.

In this section we find (asymptotic) lower bounds of $D_N(T)$ under the conditions that as $N \rightarrow \infty$,
\begin{eqnarray} \label{gam}
\log \gamma & \sim & N^\zeta \mbox{ for some } 0 < \zeta < 1, \\ \label{pNb}
p & \sim & N^{-\beta} \mbox{ for some } 0 < \beta < 1.
\end{eqnarray}
In Sections 3 and 5, 
we devise optimal detectability score stopping rules that achieve this lower bound.
In Theorem \ref{thm1} below, 
only $\beta > \frac{1-\zeta}{2}$ is considered.
For $\beta < \frac{1-\zeta}{2}$,
the detectability score stopping rules achieve asymptotic detection delay of 1,
and are hence optimal.
 
For $\frac{1-\zeta}{2} < \beta < 1-\zeta$, 
the detection delay lower bound grows logarithmically with $N$.
The proportionality constant is 
$$\rho(\beta,\zeta) = \left\{ \begin{array}{ll} \beta-\frac{1-\zeta}{2} & \mbox{ if } \frac{1-\zeta}{2} < \beta \leq \frac{3(1-\zeta)}{4}, \cr
(\sqrt{1-\zeta}-\sqrt{1-\zeta-\beta})^2 & \mbox{ if } \frac{3(1-\zeta)}{4} < \beta < 1-\zeta. \end{array} \right.
$$
This is a two-dimensional extension of the Donoho-Ingster-Jin constants $\rho(\beta) := \rho(\beta,0)$,
which has appeared in connection with sparse normal mixture detection, see \cite{DJ04, Ing97, Ing98}. 
The extension results from the additional difficulty of detecting a normal mean shift 
when there are multiple comparisons,
here for sequential change-point detection, 
and in \cite{CW15} for fixed-sample change-point detection.

\begin{thm} \label{thm1}
Let $T$ be a stopping rule such that {\rm ARL}$(T) \geq \gamma$,
with $\gamma$ satisfying {\rm (\ref{gam})}.

{\rm (a)} If {\rm (\ref{pNb})} holds with $\frac{1-\zeta}{2} < \beta < 1-\zeta$, then
\begin{equation} \label{DN1}
\liminf_{N \rightarrow \infty} \frac{D_N(T)}{\log N} \geq 2 \mu^{-2} \rho(\beta,\zeta).
\end{equation}

{\rm (b)} If {\rm (\ref{pNb})} holds with $\beta >1-\zeta$, then
\begin{equation} \label{DN2}
\liminf_{N \rightarrow \infty} \frac{\log D_N(T)}{\log N} \geq \beta+\zeta-1.
\end{equation}
\end{thm}

The phase transition between logarithmic and polynomial growth of the detection delay boundary is at $N^{1-\beta}=N^{\zeta}$,
that is, 
at $\# \cN \doteq \log \gamma$. 
By Theorem \ref{thm1}(a), 
for larger $\# \cN$ the detection delay lower bound grows at a $\log N$ rate.
By Theorem \ref{thm1}(b),
for smaller $\# \cN$ the lower bound is roughly $(\log \gamma)/\# \cN$.
The detection delay lower bound in the logarithmic domain [Theorem \ref{thm1}(a)] is closely linked to the Donoho-Ingster-Jin detection boundary for sparse normal mixture detection,
whereas the lower bound in the polynomial domain [Theorem \ref{thm1}(b)] lies in the framework of the classical lower bound established by Lorden (1971) for $N$ fixed as $\gamma \rightarrow \infty$.

We shall first establish the connection between Theorem \ref{thm1}(a) and the Donoho-Ingster-Jin detection boundary $\sqrt{2 \rho(\beta) \log N}$.
Let $t \geq \nu \geq 1$ and $k=t-\nu+1$.
If $p \sim N^{-\beta}$, 
$\frac{1}{2} < \beta < 1$,
then as 
$$k^{-1/2} \sum_{i=\nu}^t X_{ni} \sim \left\{ \begin{array}{l} \mbox{N(0,1) under } P_\infty, \cr
(1-p) \mbox{N}(0,1)+p \mbox{N}(\mu \sqrt{k},1) \mbox{ under } P_{\nu}, \cr \end{array} \right. \qquad 1 \leq n \leq N, 
$$
sparse normal mixture detection theory dictates that $k$ should satisfy
$$\mu \sqrt{k} \geq [1+o(1)] \sqrt{2 \rho(\beta) \log N} \mbox{ (i.e. } k \geq [2 \mu^{-2} \rho(\beta)+o(1)] \log N),
$$
in order for it to be possible that the sum of Type I and II error probabilities goes to zero,
when testing $P_{\nu}$ against $P_{\infty}$ with observations up to time $t$.
By (\ref{DNT}) this leads to 
\begin{equation} \label{D1}
D_N(T) \geq [2 \mu^{-2} \rho(\beta)+o(1)] \log N,
\end{equation}
for any stopping rule $T$ satisfying ARL$(T) \geq \gamma$ with $\gamma/\log N \rightarrow \infty$.
What Theorem \ref{thm1}(a) says is that under (\ref{gam}) with $\zeta$ small enough ($<1-\beta$),
$\log N$ detection is still possible with a larger asymptotic constant.

The link between Theorem 1(b) and the classical lower bound formula of Lorden is best established via the inequality in Mei \cite[Prop 2.1]{Mei06},
that for $N$ fixed,
\begin{equation} \label{R3}
D_N(T) \geq 2 \mu^{-2} \tfrac{\log \gamma}{\# \cN} + O(1) \mbox{ as } \gamma \rightarrow \infty.
\end{equation}
Theorem \ref{thm1}(b) says that for $\log \gamma \gg \# \cN$ ($\sim N^{1-\beta}$), 
the right-hand side of (\ref{R3}) gives the correct order for the attainable detection delay. 
When $\# \cN \gg \log \gamma$,
the right-hand side of (\ref{R3}) does not provide the correct order for the attainable detection delay as we have already noted in the previous paragraph situations under which a $\log N$ detection delay is required.
Therefore the $O(1)$ in (\ref{R3}) is more appropriately $O(\log N)$,
if the dependence on $N$ in $O(1)$ is made explicit.
What Theorem \ref{thm1} also says is that the transition is sharp.
Once we get out of the classical $(\log \gamma)/(\# \cN)$ domain,
we fall into the $\log N$ domain,
there are no intermediate asymptotics.

\section{Optimal detection using detectability score}

The detectability score stopping rule is motivated by the MLR stopping rules of Xie and Siegmund \cite{XS13}.
In their formulation Xie and Siegmund consider firstly the ideal situation in which $p$ and $\mu$ are known.
The most powerful test at time $t$, 
for testing the hypothesis that change-point $\nu=s$ for some $s \leq t$, 
is the log likelihood ratio
$$\ell_{\bullet st} := \sum_{n=1}^N \ell_{nst}, \mbox{ where } \ell_{nst} = \log(1-p+p e^{\mu S_{nst}-k \mu^2/2}),
$$
with $k=t-s+1$ and $S_{nst} = \sum_{i=s}^t X_{ni}$. 

Since the change-point $\nu$ is unknown,
they suggest to maximize $\ell_{\bullet st}$ over~$s$.
The unknown $\mu$ (or more precisely $\mu_n$) in $\ell_{nst}$ is substituted by $S_{nst}^+/k$, 
and a small $p_0$ is substituted for the unknown $p$. 
In summary their stopping rule can be expressed as
\begin{equation} \label{TXS}
T_{\rm XS}(p_0) = \inf \Big\{ t: \max_{k= t-s+1 \in \cK} \wht \ell_{\bullet st}(p_0) \geq b \Big\},
\end{equation}
where $\wht \ell_{\bullet st}(p_0) = \sum_{n=1}^N \wht \ell_{nst}(p_0)$ and 
$$\wht \ell_{nst}(p_0) = \log(1-p_0+p_0 e^{(Z^+_{nst})^2/2}), \quad Z_{nst} = S_{nst}/\sqrt{k}.
$$
The set $\cK$ in (\ref{TXS}) refers to a pre-determined set of window sizes.
By applying nonlinear renewal theory, 
Xie and Siegmund derive accurate analytical approximations of ARL($T$) and $D_N(T)$ 
for $T = T_{\rm XS}(p_0)$ and related stopping rules. 

Our stopping rule is also a mixture likelihood ratio but based instead on the limits of detectability.
Let  
\begin{equation} \label{TS}
T_S(p_0) = \inf \Big\{ t: \max_{k=t-s+1 \in \cK} \sum_{n=1}^N g(Z_{nst}^+) \geq b \Big\},
\end{equation}
where $g(z) = \log[1+p_0(\lambda e^{z^2/4}-1)]$ and $\lambda = 2(\sqrt{2}-1)$.
Following Lai \cite{Lai95}, 
we consider window sizes 
\begin{equation} \label{cK}
\cK = \{ 1, \ldots, k_1 \} \cup \{ \lfloor r^j k_1 \rfloor: j \geq 1 \}, \quad k_1 \geq 1, \ r > 1.
\end{equation}

\begin{thm} \label{thm2} 
Consider stopping rule $T_S(p_0)$,
$0 < p_0 \leq 1$,
with window sizes {\rm (\ref{cK})}. 
If {\rm ARL}$(T_S(p_0)) = \gamma$, 
then threshold $b \leq \log(4 \gamma^2+2 \gamma)$.
In addition,
if {\rm (\ref{gam}), (\ref{pNb})} hold and $k_1/\log N \rightarrow \infty$, 
$p_0 = c[(\log \gamma)/N]^{1/2}$ for some $c> 0$,
then the following hold as $N \rightarrow \infty$.

{\rm (a)} If $\beta < \frac{1-\zeta}{2}$, 
then $D_N(T_S(p_0)) \rightarrow 1$.

{\rm (b)} If $\frac{1-\zeta}{2} < \beta < 1-\zeta$, 
then
\begin{equation} \label{logDN}
\frac{D_N(T_S(p_0))}{\log N} \rightarrow 2 \mu^{-2} \rho(\beta,\zeta).
\end{equation}

{\rm (c)} If $\beta > 1-\zeta$, 
then
$$\frac{\log D_N(T_S(p_0))}{\log N} \rightarrow \beta+\zeta-1.
$$
\end{thm} 

{\sc Remarks}. 
Instead of (\ref{gam}),
we can model $\gamma$ growing slowly with $N$ by assuming that 
\begin{equation} \label{oNe}
\gamma/\log N \rightarrow \infty, \qquad \log \gamma = o(N^\epsilon) \mbox{ for all } \epsilon > 0.
\end{equation}
Consider the stopping rule $T_S(p_0)$ with $p_0 = cN^{-\frac{1}{2}}$ for some $c>0$.
Under (\ref{pNb}) and (\ref{oNe}), 
the asymptotic (\ref{logDN}) holds with $\zeta=0$,
and the stopping rule is optimal in view of (\ref{D1}).  

\medskip 
We shall provide some intuition here on the detectability score transformation $g$.
Consider an i.i.d sample $Z_1, \ldots, Z_N$ that is distributed as N(0,1) under the null hypothesis $H_0$. 
If $w_N \rightarrow \infty$ with $w_N= o(\sqrt{\log N})$, 
then $\# \{ n: Z_n \geq w_N \}/N$ is asymptotically normal with mean $\alpha_N$ and variance $\alpha_N/N$, 
where $\alpha_N = P_0 \{ Z_n \geq w_n \} = \int_{w_N}^\infty (2 \pi)^{-1/2} e^{-z^2/2} dz$. 

Therefore under any alternative hypothesis $H_1$, 
$\sqrt{\alpha_N/N}$ is the minimum deviation of $P_1 \{ Z_n  \geq w_N \}$ from $\alpha_N$ that is detectable.
Since $\alpha_N$ is essentially $e^{-w_N^2/2}$ (up to logarithmic terms),
the minimum detectable deviation is $e^{-w_N^2/4}/\sqrt{N}$.
That is,
a mixture of N(0,1) and a small $p_0 = cN^{-\frac{1}{2}}$ fraction of N(0,2) is at the threshold of detectability.
The detectability score transformation $g$ is essentially the likelihood ratio between the mixture with $p_0$ fraction N(0,2),
and the null distribution.
The factor $(\log \gamma)^{1/2}$ in the optimal choice of $p_0$ in the statement of Theorem \ref{thm2} adjusts for the additional difficulty of each detection due to the multiple comparison effects of large $\gamma$.

It is straightforward to check that the detectability score $\sum_{n=1}^N g(Z_{nst}^+)$ in (\ref{TS}) is indeed the log likelihood ratio for testing
$Z_{1st}^+, \ldots, Z_{Nst}^+$ i.i.d. N(0,1)$^+$
[the distribution of $Z^+$ when $Z \sim$ N(0,1)] against the alternative that $Z_{1st}^+, \ldots, Z_{Nst}^+$ are i.i.d. 
$$(1-p_0) N(0,1)^+ + p_0 [\tfrac{\lambda}{\sqrt{2}} \mbox{HN}(0,2) + (1-\tfrac{\lambda}{\sqrt{2}}) \delta_0],
$$
where $\delta_0$ denotes a point mass at zero and HN(0,2) the half-normal distribution with density $\pi^{-1/2} e^{-z^2/4}$ on $z>0$.
The value $\lambda=2(\sqrt{2}-1)$ is chosen for convenience, 
so that $g$ is continuous at 0. 
The optimality of $T_S(p_0)$ in Theorem \ref{thm2} does not require the selection of this specific $\lambda$.

\section{Numerical study}

\begin{table}[t] \label{tab1}
\begin{tabular}{c|lc|lc}
& \multicolumn{2}{|c|}{$N=100$} & \multicolumn{2}{|c}{$N=10^4$} \cr
Test & $b$ & ARL & $b$ & ARL \cr \hline
max & 12.8 & 5041 & 15.9 & 4930 \cr
Mei & 88.5 (106.8) & 4997 & 5640 (8722) & 4909 \cr 
Mei($N^{-\frac{1}{2}}$) & 3.48 (9.81) & 4994 & 3.03 (8.93) & 4973 \cr
Mei($3N^{-\frac{1}{2}}$) & 5.02 (9.61) & 4976 & 2.31 (6.97) & 5017 \cr
S($N^{-\frac{1}{2}}$) & 4.25 (18.42) & 5066 & 14.49 (18.42) & 5121 \cr
S($3N^{-\frac{1}{2}}$) & 6.30 (18.42) & 5195 & 17.21 (18.42) & 4986 

\end{tabular}
\caption{Thresholds $b$ for stopping rules calibrated to ARL $\doteq 5000$. 
The upper bounds of the thresholds,
as given in the statement of Theorems \ref{thm2} and \ref{thm3},
are in brackets.} 
\end{table}

In addition to (\ref{TXS}), 
Xie and Siegmund introduce the stopping rule  
\begin{equation} \label{TLR}
T_{\rm LR}(p_0) = \inf \Big\{ t: \max_{k=t-s+1 \in \cK} \sum_{n=1}^N (\mu_0 S_{nst} - k \mu_0^2/2 + \log p_0)^+ \geq b \Big\}.
\end{equation}
This like (\ref{TXS}) is motivated by the most powerful likelihood ratio test, 
but with $\mu$ substituted by a pre-determined $\mu_0$ rather than $S_{nst}^+/k$.  
It bears resemblance to Mei's stopping rule
\begin{equation} \label{TMei}
T_{\rm Mei} = \inf \Big\{ t: \sum_{n=1}^N \max_{0 < s \leq t} (\mu_0 S_{nst}-k \mu_0^2/2)^+ \geq b \Big\},
\end{equation}
with the important difference of an additional $\log p_0$ term in (\ref{TLR}) that suppresses the contributions of low scoring data streams.

Another key difference is that the sum lies outside the max in (\ref{TMei}) whereas in $T_{\rm LR}$ (and $T_{\rm XS}$, $T_S$), 
the sum lies inside the max. 
This confers advantage to Mei's stopping rule when the change-point $\nu$ (or $\nu_n$) differs across data streams.
We investigate this in Section 5 where we also propose an extension of Mei's stopping rule, 
denoted by $T_{\rm Mei}(p_0)$, 
that like (\ref{TLR}) weighs down the contributions from non-signal data streams.

In our numerical study, 
we benchmark the detectability score stopping rule against the above stopping rules and the max rule 
\begin{equation} \label{Tmax}
T_{\rm max} = \inf \{ t: \max_{0 < s \leq t} \max_{1 \leq n \leq N} (Z_{nst}^+)^2/2 \geq b \}.
\end{equation}
As in \cite{XS13}, 
we select $N=100$, 
$\mu=1$ and $\# \cN$ ranging from 1 to 100.
The thresholds $b$ are calibrated to average run length 5000.
The set of window sizes chosen is $\cK = \{ 1,\ldots, 200 \}$, 
and for Mei's stopping rule and $T_{\rm LR}$ we select $\mu_0=1$.

We consider $p_0=0.1(=N^{-\frac{1}{2}})$ for the detectability score stopping rule $T_S$, 
corresponding to the optimal choice under (\ref{oNe}). 
Another selection is $p_0 = 0.3 \{ \doteq [(\log \gamma)/N]^{1/2} \}$,
which is optimal under (\ref{gam}). 
It is interesting that in \cite{XS13}, 
the ``optimal" $p_0=N^{-\frac{1}{2}}$ is chosen for $T_{\rm XS}$ and $T_{\rm LR}$ in the numerical study.  

We conduct 500 Monte Carlo trials for the estimation of each average run length and detection delay.
The thresholds for the stopping rules are in Table 1, 
the detection delays in Table 2.
In Table 2 the simulation outcomes below the horizontal line are new,  
the outcomes above are reproduced from \cite[Table 5]{XS13}.  

\begin{table}[t] \label{tab2}
\begin{tabular}{c|ccccccc} 
& \multicolumn{7}{c}{$\# \cN$} \cr
Test & 1 & 3 & 5 & 10 & 30 & 50 & 100 \cr \hline
max & 25.5 & 18.1 & 15.5 & 12.6 & 9.6 & 8.6 & 7.2 \cr
XS(1) & 52.3 & 18.7 & 12.2 & 6.7 & 3.0 & 2.3 & 2.0 \cr
XS(0.1) & 31.6 & 14.2 & 10.4 & 6.7 & 3.5 & 2.8 & 2.0 \cr
LR(0.1) & 29.1 & 13.4 & 9.8 & 7.1 & 4.6 & 4.0 & 3.4 \cr
LR(1) & 82.0 & 27.2 & 15.5 & 6.8 & 3.0 & 2.3 & 2.0 \cr 
Mei & 53.2 & 23.0 & 15.7 & 9.6 & 4.9 & 3.8 & 3.0 \cr \hline
Mei(0.1) & 26.4 & 14.6 & 10.8 & 7.7 & 4.5 & 3.4 & 2.3 \cr
Mei(0.3) & 34.3 & 15.9 & 11.8 & 7.6 & 4.1 & 3.1 & 2.0 \cr 
$S$(0.1) & 26.8 & 13.4 & 9.6 & 6.4 & 2.8 & 2.0 & 1.1 \cr
$S$(0.3) & 32.6 & 14.0 & 9.5 & 5.6 & 2.3 & 1.5 & 1.0 \cr \hline
s.e. & 0.9 & 0.3 & 0.1 & 0.1 & 0.1 & 0.1 & 0.1 
\end{tabular}
\caption{Detection delays when $\# \cN$ $($out of $N=100)$ data streams undergo distribution changes.
Entries in the last row are standard error upper bounds. }
\end{table}

We see that with a few understandable exceptions, 
the detectability score stopping rules $T_S(0.1)$ and $T_S(0.3)$ have smaller detection delays compared to their competitors over the full range of $\# \cN$.
This justifies the application of the detectability score stopping rules for a relatively small $N=100$. 

Following the recommendation of a referee,
we conduct a second numerical exercise for a larger $N=10^4$,
with $\# \cN$ ranging from 1 to $10^4$.
As in the earlier simulation study,
we select $\mu=\mu_0=1$, 
ARL $=5000$ and ${\cal K}= \{1, \ldots, 200 \}$.
The detection thresholds are in Table 1,
the detection delays in Table 3.
We see again that except for $\# \cN=1$ when $T_{\max}$ is superior,
the detection score stopping rules $T_S(p_0)$ for $p_0 = 0.01 (=N^{-\frac{1}{2}})$ and 0.03 $\{ \doteq [(\log \gamma)/N]^{1/2} \}$ have the smallest detection delays.

\section{Detectability of Mei's stopping rule}

\begin{table}[t] \label{tab2a}
\begin{tabular}{c|ccccc} 
& \multicolumn{5}{c}{$\# \cN$} \cr
Test & 1 & 10 & $10^2$ & $10^3$ & $10^4$ \cr \hline
max & 32.7 & 18.6 & 13.9 & 11.1 & 9.4 \cr
Mei & 246.5 & 46.7 & 12.0 & 4.0 & 1.0 \cr
Mei(0.01) & 39.7 & 16.7 & 8.8 & 4.0 & 2.0 \cr
Mei(0.03) & 53.7 & 18.6 & 9.0 & 4.0 & 2.0 \cr 
$S$(0.01) & 37.7 & 13.3 & 4.5 & 1.0 & 1.0 \cr
$S$(0.03) & 49.3 & 13.7 & 3.9 & 1.0 & 1.0 \cr
\hline
s.e. & 4.0 & 0.3 & 0.1 & 0.1 & 0.1  
\end{tabular}
\caption{Detection delays when $\# \cN$ (out of $N=10^4$) data streams undergo distribution changes.
Entries in the last row are standard error upper bounds. }
\end{table}

As mentioned earlier there is no implicit assumption that the distribution changes occur simultaneously when applying Mei's stopping rule (\ref{TMei}).
Another advantage is the efficient recursive computation of the stopping rule.
However this recursive computation comes with the price of information loss.  
In this section we improve Mei's stopping rule by applying a detectability score transformation on each CUSUM score.
Due to the information loss,
optimality is possible only for specific $\mu_0$.

Let $R_{nt}$ be the CUSUM score of the $n$th detector at time $t$, satisfying  
\begin{equation} \label{Rnt}
R_{n0}=0, \quad R_{nt} = (R_{n,t-1}+\mu_0 X_{nt}-\mu_0^2/2)^+, \quad t \geq 1.
\end{equation}
Define 
\begin{equation} \label{TMp}
T_{\rm Mei}(p_0) = \inf \Big\{ t: \sum_{n=1}^N g_M(R_{nt}) \geq b  \Big\},
\end{equation}
with the detectability score transformation
\begin{equation} \label{gMx}
g_M(x) = \log [1+p_0(\lambda_M e^{x/2}-1)], \quad \lambda_M > 0.
\end{equation}
This is an extension of Mei's test, 
for $T_{\rm Mei}(1)$ is equivalent to $T_{\rm Mei}$.
Let $\xi = \lim_{t \rightarrow \infty} E_\infty e^{R_{nt}/2}$ and define
$$D_{N,k}(T) = \sup_{k \leq \nu < \infty} E_\nu(T-\nu+1|T \geq \nu).
$$

\begin{thm} \label{thm3}
Consider stopping rule $T_{\rm Mei}(p_0)$,
$0 < p_0 \leq 1$. 
Let $u=\log[1+p_0(\lambda_M \xi-1)]$.
If {\rm ARL}$(T_{\rm Mei}(p_0)) = \gamma$,
then threshold $b \leq Nu+\log(4 \gamma)$.
In addition,
if {\rm (\ref{gam}), (\ref{pNb})} hold and $p_0 = c[(\log \gamma)/N]^{1/2}$ for some $c> 0$,
then the following hold as $N \rightarrow \infty$.  

{\rm (a)} If $\frac{1-\zeta}{2} < \beta \leq \frac{3(1-\zeta)}{4}$ and $\mu_0 = 2 \mu$, 
then 
\begin{equation} \label{3a}
\frac{D_{N,K_N}(T_{\rm Mei}(p_0))}{\log N} \rightarrow 2 \mu^{-2} \rho(\beta,\zeta),
\end{equation}
for $K_N = 2 \mu^{-2} (1-\zeta-\beta) \log N$.

\smallskip
{\rm (b)} If $\frac{3(1-\zeta)}{4} < \beta < 1-\zeta$ and $\mu_0 = \mu \sqrt{\frac{1-\zeta}{\rho(\beta,\zeta)}}$, 
then {\rm (\ref{3a})} holds for $K_N = 2 \mu^{-2} \rho(\beta,\zeta) \log N$. 
\end{thm}

{\sc Remarks}. 
1. In Theorem \ref{thm3} ``optimality" occurring when $\mu_0 > \mu$ is a consequence of a small subset of $\cN$ dominating the score contributions, 
after the detectability score transformations have been applied.

2. Notice the weaker (\ref{3a}) instead of (\ref{logDN}). 
The extra initial delay is needed for the CUSUM scores $R_{nT}$ for $n \not\in \cN$ to reach their stationary values and not pull down the total score. 
In that sense the detection delay criterion may be disadvantageous to the extended Mei's stopping rule (and hence Mei's test stopping rule itself) since in practice we seldom expect the change-point $\nu$ to be that close to 0.

\medskip
To highlight the unique characteristics of the extended Mei's stopping rule (\ref{TMp}) in dealing with staggered change-points, 
we conduct a numerical study with $\mu_{nt} = \mu {\bf I}_{ \{ t \geq n \}}$ in place of (\ref{munt}).
That is the $n$th data stream undergoes a distribution change at time $n$. 
As in Section 4 the stopping rules are calibrated to average run length of 5000, 
for $N=100$ detectors, 
and with $\mu_0=1$. 
The thresholds $b$ for $T_{\rm Mei}(p_0)$ are in Table 1 (Section 4),
the detection delays in Section 4.
We select $\lambda_M = 0.64$,
this will be explained later. 
By detection delay we shall mean the expected stopping time when $\mu_{nt} = \mu {\bf I}_{\{ t \geq n \}}$.

\begin{table}[t] \label{tab3}
\begin{tabular}{c|ccccc}
$\mu$ & Mei  & Mei(0.1) & Mei(0.3) & $S(0.1)$ & $S(0.3)$ \cr \hline
0.5 & 20.7 & 21.2 & 20.6 & 23.0 & 20.7 \cr
0.7 & 15.5 & 15.4 & 15.1 & 16.0 & 14.9 \cr
1.0 & 11.9 & 10.9 & 11.1 & 10.9 & 10.4 \cr
1.3 & 10.0 & 8.7 & 9.0 & 8.0 & 7.9
\end{tabular}
\caption{Detection delays for staggered distribution changes.
The standard errors are not more than 0.2.} 
\end{table}

We see from Tables 2 (Section 4) and 4 that $T_{\rm Mei}(0.1)$ and $T_{\rm Mei}(0.3)$ have smaller detection delays compared to $T_{\rm Mei}$, 
almost uniformly over $\# \cN$ and $\mu$.
In Table 3 (for $N=10^4$),
$T_{\rm Mei}(0.01)$ and $T_{\rm Mei}(0.03)$ are superior to $T_{\rm Mei}$ for $\# \cN \leq 100$.
Hence applying detectability score transformations on the CUSUM scores improves Mei's stopping rule in general,
the noise suppression on data streams that do not undergo distribution change is indeed effective.
In Table 4 we see that in general $T_S(p_0)$ performs better than $T_{\rm Mei}(p_0)$ when $\mu \geq 1$ but the reverse is true when $\mu < 1$.
This is consistent with the prediction in Theorem 3 of $T_{\rm Mei}(p_0)$ performing better for $\mu < \mu_0$.  

We end this section with explanations of the choice of the detectability score transformation (\ref{gMx}) and choice of $\lambda_M$.
It follows from renewal theory,
see for example Siegmund \cite[eq8.49]{Sie85}, 
that 
\begin{equation} \label{limR} 
\lim_{t \rightarrow \infty} P_\infty \{ R_{nt} \geq x \} \sim \alpha e^{-x} \mbox{ as } x \rightarrow \infty,
\end{equation}
for $\alpha = 2 \mu_0^{-2} \exp[-2 \sum_{j=1}^\infty j^{-1} \Phi(-\mu_0 \sqrt{j}/2)]$.
Therefore the tails of $R_{nt}$ under $P_\infty$ are like that of an i.i.d. sample from $G_1 := (1-\alpha) \delta_0 + \alpha {\rm Exp}(1)$,
where $\delta_0$ denotes a point mass at 0 and Exp($\theta$) the exponential distribution with mean $\theta$.

For large $x$ (smaller than $\log N$) and $t$,
$\# \{ n: R_{nt} \geq x \}/N$ is asymptotically normal with mean $\alpha e^{-x}$ and variance $\alpha e^{-x}/N$.
Hence the minimum detectable difference of $P \{ R_{nt} \geq x \}$ is $e^{-x/2}/\sqrt{N}$.
The distribution at the limit of detectability is therefore $G^* := (1-p_0)G_1+p_0 G_2$,
where $G_2=(1-\omega) \delta_0 + \omega {\rm Exp}(2)$ for some $0 < \omega < 1$, 
and $p_0$ is of order $N^{-\frac{1}{2}}$.
The detectability score transformation $g_M$ [see (\ref{gMx})],
with $\lambda_M=\frac{1}{1+\alpha} (=0.64$ for $\mu_0=1$), 
is the log likelihood ratio between $G^*$ and $G_1$,
with $\omega$ selected so that $g_M$ is continuous at 0.
We emphasize however that this is for convenience, 
optimality in Theorem \ref{thm3} is not restricted to this choice of $\lambda_M$.

\section{Proof of Theorem \ref{thm1}}

To help the reader,
we summarize below the definitions of the probability measures used in the proofs of Theorems \ref{thm1}--\ref{thm3} in this and the next two sections.

\begin{enumerate}
\item $P_s$ ($E_s$): 
This is the probability measure (expectation) under which an arbitrarily chosen data stream has probability $(1-p)$ that all observations are (i.i.d.) N(0,1),
and probability $p$ that observations are N(0,1) before time $s$,
N($\mu$,1) at and after time $s$.
In particular,
if 
\begin{enumerate}
\item $s=\infty$,
then with probability 1 all observations are N(0,1).

\item $s=1$, 
then an arbitrarily chosen data stream has probability $(1-p)$ that all observations are N(0,1),
and probability $p$ that all observations are N($\mu$,1).
\end{enumerate}

\item $P$ ($E$): 
This is the probability measure (expectation) under which $Y, Y_1, Y_2, \ldots$ are i.i.d. N(0,1) random variables. 
\end{enumerate}

We preface the proof of Theorem \ref{thm1} with the following lemmas.
Lemma \ref{lem1} is well-known,
see for example (3.3) of Lai \cite{Lai95}.

\begin{lem} \label{lem1}
Let $k \geq 1$.
If $T$ is a stopping rule such that $E_\infty T \geq \gamma$, 
then $P_\infty \{ T \geq s+k|T \geq s \} \geq 1-k/\gamma$ for some $s \geq 1$.
\end{lem}

\medskip
Recall the sum $S_{nst} = \sum_{i=s}^t X_{ni}$ and the log likelihood ratio
$$\ell_{\bullet st} = \sum_{n=1}^N \ell_{nst}, 
\mbox{ where } \ell_{nst} = \log (1-p+pe^{\mu S_{nst}-k \mu^2/2}), \quad k=t-s+1.
$$

\begin{lem} \label{lem2}
If we can find $b(=b_N)$ and $k(=k_N)$ such that  
\begin{eqnarray} \label{const}
& & P_\infty \{ \ell_{\bullet 1k} \geq b \} (=P_\infty \{ \ell_{\bullet st} \geq b \}) \geq k/\gamma, \\
\label{(2)}
& & P_1 \{ \ell_{\bullet 1k} \geq b \} (=P_s \{ \ell_{\bullet st} \geq b \}) \rightarrow 0,
\end{eqnarray}
then $D_N(T) \geq [1+o(1)]k$ for any stopping rule $T$ satisfying $E_\infty T \geq \gamma$.
\end{lem}

{\sc Proof}. 
Let $T$ satisfies $E_\infty T \geq \gamma$, 
and let $b, k$ satisfy (\ref{const}) and (\ref{(2)}).
By Lemma \ref{lem1} we can find $s$ satisfying 
\begin{equation} \label{7*}
P_\infty \{ T \geq s+k | T \geq s \} \geq 1-k/\gamma.
\end{equation}
Let $P_\infty^* \{ \cdot \} = P_\infty \{ \cdot | T \geq s \}$ and $P_s^* \{ \cdot \} = P_s \{ \cdot | T \geq s \}$.

Let $t=s+k-1$, 
and consider the test, conditioned on $T \geq s$,
of 
$$\begin{array}{rl}
H_0: & X_{nu} \sim \mbox{ N(0,1) for } 1 \leq n \leq N, 1 \leq u \leq t, \cr
\mbox{vs } H_s: & X_{nu} \sim \mbox{N}(\mu {\bf I}_{\{ u \geq s, n \in \cN \}},1) \mbox{ for  } 1 \leq n \leq N, 1 \leq u \leq t, \cr
& \mbox{ with } {\bf I}_{\{ n \in \cN \}} \sim \mbox{ Bernoulli}(p).
\end{array}
$$
By (\ref{7*}) the test ``reject $H_0$ if $T<s+k$, 
accept $H_0$ otherwise'' has Type I error probability not exceeding $k/\gamma$.
By (\ref{const}) the likelihood ratio test rejecting $H_0$ when $\ell_{\bullet st}$ exceeds $b$ has Type I error probability at least $k/\gamma$,
and hence by the Neyman-Pearson Lemma,
it is at least as powerful as the test based on $T$.
That is
\begin{equation} \label{Psstar}
P_s^* \{ \ell_{\bullet st} \geq b \} \geq P_s^* \{ T < s+k \}.
\end{equation}
A key observation here is that the conditioning on $\{ T \geq s \}$ does not affect the distribution of $X_{nu}$ for $u \geq s$ under either $H_0$ or $H_s$.
Therefore by (\ref{Psstar}),
\begin{eqnarray*} 
D_N(T) & \geq & E_s(T-s+1|T \geq s) \geq k P_s^* \{ T \geq s+k \} \cr
& \geq & k P_s^* \{ \ell_{\bullet st} < b \} = k P_s \{ \ell_{\bullet st} < b \},
\end{eqnarray*}
and we conclude $D_N(T) \geq [1+o(1)]k$ from (\ref{(2)}). 
$\wbox$

\begin{lem} \label{lem2.5}
If $k$ is such that $\log k=o(N^\zeta)$ and
\begin{equation} \label{p0}
P_1 \{ \ell_{\bullet 1k} \geq 2N^\zeta/3 \} \rightarrow 0,  
\end{equation}
then {\rm (\ref{const})} and {\rm (\ref{(2)})} follow from selecting $b$ satisfying
\begin{equation} \label{p3}
P_1 \{ 2N^\zeta/3 \geq \ell_{\bullet 1k} \geq b \} = \exp(-N^\zeta/4).
\end{equation}
\end{lem}

{\sc Proof}. 
It follows from (\ref{p0}) and (\ref{p3}) that (\ref{(2)}) holds.
Moreover since $\ell_{\bullet 1k}$ is the log change of measure between $P_1$ and $P_\infty$ at time $k$,
\begin{eqnarray*}
& & P_\infty \{ \ell_{\bullet 1k} \geq b \} \geq P_\infty \{ 2N^\zeta/3 \geq \ell_{\bullet 1k} \geq b \} \cr
& = & E_1(e^{-\ell_{\bullet 1k}} {\bf I}_{\{ 2N^\zeta/3 \geq \ell_{\bullet 1k} \geq b \}}) \geq \exp(-2N^\zeta/3) P_1 \{ 2N^\zeta/3 \geq \ell_{\bullet 1k} \geq b \},
\end{eqnarray*}
and (\ref{const}) follows from (\ref{p3}) since $\log (\gamma/k) \sim N^\zeta$.
$\wbox$

\medskip
In view of Lemmas \ref{lem2} and \ref{lem2.5},
to prove Theorem \ref{thm1} it suffices to check (\ref{p0}) for 
\begin{equation} \label{k=}
k = \left\{ \begin{array}{ll} 
\lfloor (1-\delta) 2 \mu^{-2} \rho(\beta,\zeta) \log N \rfloor & \mbox{ if } \frac{1-\zeta}{2} < \beta < 1-\zeta, \cr 
\lfloor \delta N^{\beta+\zeta-1} \rfloor & \mbox{ if } \beta > 1-\zeta, \end{array} \right.
\end{equation}
with $\delta > 0$ small. 
Motivations behind the above choices of $k$ are given in Appendix A.

Let $Z_{nk} = S_{n1k}/\sqrt{k}$ and 
\begin{equation} \label{ellnk}
\ell_{nk} (=\ell_{n1k}) = \log(1-p+p e^{Z_{nk} \mu \sqrt{k}-k \mu^2/2}).
\end{equation}
Note that $Z_{nk}$,
$1 \leq n \leq N$, 
are i.i.d. N(0,1) under $P_\infty$, 
and i.i.d. $(1-p)$N(0,1)$+p$N($\mu \sqrt{k}$,1) under $P_1$.
More specifically, 
$Z_{nk}$ has the distribution of $Y \sim$ N(0,1) if $n \not\in \cN$,
and the distribution of $Y+\mu \sqrt{k}$ if $n \in \cN$.
Hence conditioned on $n \not\in \cN$,
$\ell_{nk}$ has the distribution of
\begin{equation} \label{ell0}
\ell_0 = \log(1-p+pe^{Y \mu \sqrt{k}-k \mu^2/2}),
\end{equation}
whereas conditioned on $n \in \cN$,
$\ell_{nk}$ has the distribution of
\begin{equation} \label{ell1}
\ell_1 = \log(1-p+pe^{Y \mu \sqrt{k}+k \mu^2/2}).
\end{equation}

\medskip
Case 1: $\frac{1-\zeta}{2} < \beta < 1-\zeta$. 
Let $\wtd \ell_{nk} = \ell_{nk} {\bf I}_{\{ Z_{nk} \leq \omega_N \}}$, 
where 
$$\omega_N = \sqrt{2(1-\zeta) \log N+2 \log \log N}.
$$
We shall check on two sub-cases that
\begin{eqnarray} \label{62a2}
\wtd \mu := E_1 \wtd \ell_{nk} & = & o(N^{\zeta-1}), \\ \label{62a3}
\sup_{1 \leq n \leq N} \wtd \ell_{nk}^+ & = & O(1), \\ \label{vdef}
E_1 \wtd \ell_{nk}^2 & = & o(N^{\zeta-1}), \\ \label{62b1}
P_1 \{ Z_{nk} > \omega_N \} & = & o(N^{\zeta-1}/\log N).
\end{eqnarray}
Note that by (\ref{62b1}) and $\max_{1 \leq n \leq N} Z_{nk} = O_p(\sqrt{\log N})$,
\begin{equation} \label{large}
\sum_{n=1}^N \ell_{nk} {\bf I}_{\{ Z_{nk} > \omega_N \}} =o_p(N^{\zeta}/\sqrt{\log N}).
\end{equation}

Recall that $\ell_{\bullet 1k} = \sum_{n=1}^N \ell_{nk}$ and let $\wtd \ell_{\bullet 1k} = \sum_{n=1}^N \wtd \ell_{nk}$.
By Chebyshev's inequality and (\ref{vdef}),
\begin{eqnarray} \label{p1}
& & P_1 \{ \wtd \ell_{\bullet 1k} - N \wtd \mu \geq N^{\zeta/2} \} \leq N^{-\zeta} E_1(\wtd \ell_{\bullet 1k}-N \wtd \mu)^2 \\ \nonumber
& = & N^{-\zeta+1} E_1 (\wtd \ell_{nk}-\wtd \mu)^2 \leq N^{-\zeta+1} E_1 \wtd \ell_{nk}^2 \rightarrow 0.
\end{eqnarray}
By (\ref{large}),
noting that $\ell_{\bullet 1k} - \wtd \ell_{\bullet 1k} = \sum_{n=1}^N \ell_{nk} {\bf I}_{\{ Z_{nk} > \omega_N \}}$,
\begin{equation} \label{pp2}
P_1 \{ \ell_{\bullet 1k} - \wtd \ell_{\bullet 1k} \geq N^\zeta/\sqrt{\log N} \} \rightarrow 0.
\end{equation}
It follows from (\ref{p1}) and (\ref{pp2}) that $P_1 \{ \ell_{\bullet 1k} \geq \wht b \} \rightarrow 0$ for $\wht b = N \wtd \mu+N^{\zeta/2}+N^\zeta/\sqrt{\log N}[=o(N^\zeta)$ by (\ref{62a2})],
hence (\ref{p0}) holds.  

\medskip
Checking (\ref{62a2})--(\ref{62b1}):

\smallskip
(a) $\frac{1-\zeta}{2} < \beta \leq \frac{3(1-\zeta)}{4}$ and $\rho(\beta,\zeta) = \beta-\frac{1-\zeta}{2}$.
By Jensen's inequality, 
$E \ell_0 \leq \log Ee^{\ell_0} = 0$,
therefore to show (\ref{62a2}), 
it suffices to show that
\begin{equation} \label{3u}
p E \ell_1^+ = o(N^{\zeta-1}).
\end{equation}
Indeed as $\log (1+x) \leq x$,
by (\ref{k=}),
\begin{equation} \label{pp}
p E \ell_1^+ \leq p^2 E e^{Y \mu \sqrt{k}+k \mu^2/2} = p^2 e^{k \mu^2} = O(N^{-2\beta+(1-\delta)(2 \beta-1+\zeta)}),
\end{equation}
and (\ref{3u}) holds.

To show (\ref{62a3}), 
note that
\begin{eqnarray} \label{63a1}
\sup_{1 \leq n \leq N} \wtd \ell_{nk}^+ \leq p e^{\omega_N \mu \sqrt{k}-\mu k^2/2} & = & pe^{\omega_N^2/2-(\omega_N-\mu \sqrt{k})^2/2} \\ \nonumber
& \sim & N^{-\beta+1-\zeta} e^{-(w_N-\mu \sqrt{k})^2/2} \log N.   
\end{eqnarray}
Express $\beta=\frac{1-\zeta}{2}+\alpha(1-\zeta)$ for some $0 < \alpha < \frac{1}{4}$.
Since $\rho(\beta,\zeta)=\alpha(1-\zeta)$ and $w_N \geq \sqrt{2(1-\zeta) \log N}$,
by (\ref{k=}) there exists $\epsilon > 0$ small such that
\begin{eqnarray} \label{wsq}
\tfrac{(w_N-\mu \sqrt{k})^2}{2 \log N} & \geq & (1-\zeta)(1-\sqrt{\alpha})^2+\epsilon \\ \nonumber
& = & (1-\zeta)[\tfrac{(1-2 \sqrt{\alpha})^2}{2}+\tfrac{1}{2}-\alpha]+\epsilon \\ \nonumber
& \geq & (1-\zeta)(\tfrac{1}{2}-\alpha)+\epsilon \\ \nonumber
& = & 1-\zeta-\beta+\epsilon.
\end{eqnarray}
Substituting (\ref{wsq}) into (\ref{63a1}) shows (\ref{62a3}).

To show (\ref{vdef}),
note that by (\ref{pp}),
\begin{equation} \label{E02}
E \ell_0^2 = O(p^2 e^{2Y \mu \sqrt{k}-k \mu^2}) = O(p^2 e^{k \mu^2}) = o(N^{\zeta-1}).
\end{equation}
Since $\beta > \frac{1-\zeta}{2}$,
\begin{equation} \label{E1n}
(\wtd \ell_{nk}^-)^2 = O(p^2) = o(N^{\zeta-1}),
\end{equation}
and (\ref{vdef}) follows from (\ref{62a3}), (\ref{3u}) and (\ref{E02}). 

Finally to show (\ref{62b1}), 
note that $P \{ Y > \omega_N \} = o(N^{\zeta-1}/\log N)$, 
and that by (\ref{wsq}),
\begin{eqnarray} \label{Pstar}
p P \{ Y+\mu \sqrt{k} > \omega_N \} & = & O(N^{-\beta} e^{-(\omega_N - \mu \sqrt{k})^2/2}) \\ \nonumber
& = & o(N^{\zeta-1}/\log N). 
\end{eqnarray}

(b) $\frac{3(1-\zeta)}{4} < \beta < 1-\zeta$ and $\rho(x,y)=(x-y)^2$, 
where $x=\sqrt{1-\zeta}$, $y=\sqrt{1-\zeta-\beta}$.
By $\log(1+v) \leq v$,
\begin{eqnarray} \label{pEl}
p E (\ell_1^+ {\bf I}_{ \{ Y+\mu \sqrt{k} \leq \omega_N \}}) & \leq &  
p^2 \int_{-\infty}^{\omega_N-\mu \sqrt{k}} \tfrac{1}{\sqrt{2 \pi}} e^{-z^2/2+z \mu \sqrt{k}+k \mu^2/2} dz \\ \nonumber
& = & p^2 e^{k \mu^2} \Phi(\omega_N-2 \mu \sqrt{k}).
\end{eqnarray}
Since $\omega_N \sim  x \sqrt{2 \log N}$, $\mu \sqrt{k}=(1-\delta)(x-y) \sqrt{2 \log N}+O(1)$ and $x > 2y$,
it follows that $\omega_N < 2 \mu \sqrt{k}$ for $\delta > 0$ small, 
and therefore
\begin{eqnarray} \label{p2}
p^2 e^{k \mu^2} \Phi(\omega_N - 2 \mu \sqrt{k}) & = & O(p^2 e^{k \mu^2-(\omega_N-2 \mu \sqrt{k})^2/2}) \\ \nonumber
& = & O(p^2 e^{\omega_N^2/2-(\omega_N-\mu \sqrt{k})^2}) \\ \nonumber
& = & O(N^{-2 \beta+x^2-2y^2-\epsilon}) = O(N^{\zeta-1-\epsilon})
\end{eqnarray}
for some $\epsilon > 0$,
since $-2\beta+x^2-2y^2 = \zeta-1$,
And since $E \ell_0 \leq 0$, (\ref{62a2}) follows from (\ref{pEl}) and (\ref{p2}). 

By the first line of (\ref{63a1}), 
$$\wtd \ell_{nk} \leq p e^{\omega_N \mu \sqrt{k}-\mu k^2/2} =O(N^{-\beta+x^2-y^2-\epsilon}) 
$$ 
for some $\epsilon > 0$, 
therefore (\ref{62a3}) holds.

Note that by (\ref{ell0}) and $\log(1+v) \leq v$,
\begin{eqnarray*}
E(\ell_0^2 {\bf I}_{\{ Y \leq \omega_N \}}) & = & O \Big( p^2 \int_{-\infty}^{\omega_N} \tfrac{1}{\sqrt{2 \pi}} e^{-z^2/2+2z \mu \sqrt{k}-k \mu^2} dz \Big) \cr
& = & O(p^2 e^{k \mu^2} \Phi(\omega_N-2 \mu \sqrt{k})),
\end{eqnarray*}
and (\ref{vdef}) follows from (\ref{62a3}), (\ref{E1n}), (\ref{pEl}) and (\ref{p2}).
It is easy to check that (\ref{Pstar}), 
and hence (\ref{62b1}), 
holds in this sub-case. 

\medskip 
Case 2: $\beta > 1-\zeta$ and $k= \lfloor \delta N^{\beta+\zeta-1} \rfloor$. 
Let
$$\wht \ell_{nk} = \left\{ \begin{array}{ll} \ell_{nk} {\bf I}_{\{ Z_{nk} \leq \sqrt{2 \log N} \}} & \mbox{ if } n \not\in \cN, \cr
\ell_{nk} & \mbox{ if } n \in \cN. \end{array} \right.
$$
Let $\wht \ell_{\bullet 1k} = \sum_{n=1}^N \wht \ell_{nk}$.
In place of (\ref{62a2}) and (\ref{vdef}), 
we shall check that for $\delta > 0$ small and $N$ large,
\begin{eqnarray} \label{52a}
E_1 \wht \ell_{nk} & \leq & N^{\zeta-1}/2, \\ \label{wl2}
E_1 \wht \ell^2_{nk} & = & o(N^{2 \zeta-1}).
\end{eqnarray}
Note that in place of (\ref{large}), 
we have
\begin{equation} \label{P1Z}
P_1 \{ Z_{nk} > \sqrt{2 \log N} \mbox{ for some } n \not\in \cN \} (= P_1 \{ \ell_{\bullet 1k} > \wht \ell_{\bullet 1k} \}) \rightarrow 0.
\end{equation}
It follows from (\ref{52a}), (\ref{wl2}) and Chebyshev's inequality,
see the arguments in (\ref{p1}), 
that $P_1 \{ \wht \ell_{\bullet 1k} \geq 2N^\zeta/3 \} \rightarrow 0$, 
hence (\ref{p0}) follows from (\ref{P1Z}). 

Check that
\begin{equation} \label{log1}
\log(1+e^x) \leq \log 2 + x^+, 
\end{equation}
and apply it on (\ref{ell1}) to show that 
$$p E \ell_1 \leq p [\log 2 + E(\log p+Y \mu \sqrt{k} + k \mu^2/2)^+] \sim \delta \mu^2 N^{\zeta-1}/2.
$$
Since $E \ell_0 \leq 0$,
(\ref{52a}) holds when $\delta < \mu^{-2}$.

Since $\sup_{n \not\in \cN} |\wht \ell_{nk}| \sim p$,
by (\ref{log1}), 
\begin{eqnarray*}
E_1 \wht \ell_{nk}^2 & \leq & p E(\log p+Y \mu \sqrt{k}+k \mu^2/2)^2 +O(p^2) \cr
& = &  O(N^{\beta+2 \zeta-2}) + O(N^{-2 \beta}),
\end{eqnarray*}
and (\ref{wl2}) holds because $\beta<1$. 

\section{Proof of Theorem \ref{thm2}}

The following lemma provides an upper bound for the threshold of the detectability score stopping rule.

\begin{lem} \label{lem3}
Consider stopping-rule $T_S(p_0)$, 
$0 < p_0 \leq 1$,
with arbitrary window-sizes ${\cal K}$.
If $b = \log(4 \gamma^2+2 \gamma)$, 
then $E_\infty T_S(p_0) \geq \gamma$.
\end{lem}

{\sc Proof}. It suffices to show that 
\begin{equation} \label{7.1}
P_\infty \{ T_S(p_0) < 2 \gamma \} \leq \tfrac{1}{2}.
\end{equation}

Let $Z_{nst} = S_{nst}/\sqrt{k}$,
$k=t-s+1$.
Since $V_{st} := \sum_{n=1}^N g(Z_{nst}^+)$ is a log likelihood ratio against $Z_{1st}, \ldots, Z_{Nst}$ i.i.d. N(0,1),
it follows from a change of measure argument that
$$P_{\infty} \{ V_{st} \geq b \} \leq e^{-b} = (4 \gamma^2+2 \gamma)^{-1}. 
$$
By Bonferroni's inequality,
$$P_{\infty} \{ T_S < 2 \gamma \} \leq \sum_{(s,t): 1 \leq s \leq t < 2 \gamma} P_{\infty} \{ V_{st} \geq b \} \leq {\lfloor 2 \gamma+1 \rfloor \choose 2} (4 \gamma^2+2 \gamma)^{-1},
$$
and (\ref{7.1}) follows. 
$\wbox$

\medskip
Assume (\ref{gam}),
(\ref{pNb}) and let $\eta = \min_{m \in J_N} P_1 \{ \sum_{n=1}^N g(Z_{n1k}^+) \geq b | \# \cN=m \}$,
where
\begin{equation} \label{JN}
J_N = \{ m: |m-Np| \leq N^{(\zeta+1)/2} \}.
\end{equation}
By the Chernoff-Hoeffding's inequality,
\begin{equation} \label{CH}
P_1 \{ \# \cN \not\in J_N \} \leq \exp(-2N^{\zeta}) = o(\gamma^{-1}).
\end{equation}
We shall show in various cases below that $\eta \rightarrow 1$ when
\begin{equation} \label{k}
k = \left\{ \begin{array}{ll} 1 & \mbox{ if } \beta < \frac{1-\zeta}{2}, \cr
\lfloor (1+\delta) 2 \mu^{-2} \rho(\beta,\zeta) \log N \rfloor & \mbox{ if } \frac{1-\zeta}{2} < \beta < 1-\zeta, \cr
M N^{\beta+\zeta-1} & \mbox{ if } \beta > 1-\zeta, \end{array} \right.
\end{equation}
for all $\delta > 0$, 
and $M$ large.
For $j \geq 1$ and $m \in J_N$, 
$P_1 \{ T_S(p_0) \geq  jk+1 | \# \cN = m \} \leq (1-\eta)^j$.
Hence by (\ref{CH}), 
\begin{equation} \label{DNTS}
D_N(T_S(p_0)) \leq k \sum_{j=0}^\infty (1-\eta)^j  + \gamma P_1 \{ \# \cN \not\in J_N \} \sim k,
\end{equation}
and the proof of Theorem \ref{thm2} is complete.

Let $V_N = \sum_{n=1}^N g(Y_n^+)$ for $Y_1, \ldots, Y_N$ i.i.d. {\rm N(0,1)}.

\begin{lem} \label{lem4}
If $p_0 \sim cN^{(\zeta-1)/2}$, 
then $P \{  V_N \geq -N^{\zeta} \} \rightarrow 1$.
\end{lem} 

{\sc Proof}. 
Let $\bar \Phi(z) = \int_z^\infty \phi(y) dy$ where $\phi(y) = \frac{1}{\sqrt{2 \pi}} e^{-y^2/2}$, 
and $\wtd g(z) = g(z) {\bf I}_{\{ z \leq w_N \}}$ where $w_N = \sqrt{2(1-\zeta) \log N - (\log \log N)/2}$.
By $\log(1+x) \sim x$ and $p_0 e^{w_N^2/4} \rightarrow 0$,
\begin{eqnarray} \label{7.5}
E \wtd g(Y_1) & \sim & c N^{(\zeta-1)/2} \int_{-\infty}^{w_N} 
(\lambda e^{z_+^2/4}-1) \phi(z) dz \\ \nonumber
& = & c N^{(\zeta-1)/2} \Big[ \tfrac{\lambda}{2}+\lambda \sqrt{2} \int_0^{w_N} \tfrac{1}{\sqrt{4 \pi}} e^{-z^2/4} dz -\Phi(w_N) \Big] \\ \nonumber
& = & c N^{(\zeta-1)/2} \{ [\tfrac{\lambda}{2}-\Phi(w_N)] +\lambda \sqrt{2} [\tfrac{1}{2}-\bar \Phi(\tfrac{w_N}{\sqrt{2}})] \}. 
\end{eqnarray}
Since $\lambda=2(\sqrt{2}-1)$ solves $\frac{\lambda}{2}+\frac{\lambda}{\sqrt{2}}=1$ and $\bar \Phi(w_N) \leq \bar \Phi(\frac{w_N}{\sqrt{2}})=o(N^{(\zeta-1)/2})$,
by (\ref{7.5}), 
\begin{equation} \label{7.6}
| E \wtd g(Y_1) | =o(N^{\zeta-1}).
\end{equation}
Since
$$E \wtd g^2(Y_1) \sim c^2 N^{\zeta-1} \int_{-\infty}^{w_N} (\lambda e^{z_+^2/4}-1)^2 \phi(z) dz = O(N^{\zeta-1} \sqrt{\log N}),
$$
and $g \geq \wtd g$,
we conclude Lemma \ref{lem4} from (\ref{7.6}) and Chebyshev's inequality. $\wbox$

\medskip Let $h(z) = g((z+\mu \sqrt{k})^+)-g(z^+) (\geq 0)$ and $H_N = \sum_{n \in \cN} h(Y_n)$.
Then
$$\eta = \min_{m \in J_N} P \{ V_N+H_N \geq b | \# \cN=m \}.
$$
In view of Lemmas \ref{lem3} and \ref{lem4}, 
to show $\eta \rightarrow 1$ and hence (\ref{DNTS}),
it suffices to show that
\begin{equation} \label{7.8}
\min_{m \in J_N} P \{ H_N \geq 4N^{\zeta} | \# \cN = m \} \rightarrow 1.
\end{equation}
We shall check (\ref{7.8}) in three cases below.
Note that $\# \cN \in J_N$ implies $\# \cN \sim N^{1-\beta}$.
For notational simplicity,
we shall let $C$ denote a generic positive constant.

\medskip Case 0: $\beta < \frac{1-\zeta}{2}$ and $k=1$.
Since $\log(1+x) \sim x$ as $x \rightarrow 0$,
$$h(z) \sim c \lambda N^{(\zeta-1)/2} (e^{(z+\mu)^2/4}-e^{z^2/4}) \geq c \lambda N^{(\zeta-1)/2} (e^{\mu^2/4}-1),
$$
uniformly over $0 \leq z \leq 1$.
Hence by LLN,
$$H_N \geq [C+o_p(1)] N^{1-\beta+(\zeta-1)/2},
$$
and (\ref{7.8}) holds because $1-\beta+\frac{\zeta-1}{2} > \zeta$.

\medskip
Case 1: $\frac{1-\zeta}{2} < \beta < 1-\zeta$ and $k= \lfloor (1+\delta) 2 \mu^{-2} \rho(\beta,\zeta) \log N \rfloor$.
We shall show (\ref{7.8}) for the following sub-cases.

\smallskip (a) $\frac{1-\zeta}{2} < \beta \leq \frac{3(1-\zeta)}{4}$ and $\rho(\beta,\zeta)=\beta-\frac{1-\zeta}{2}$.
For $\delta>0$ small and $z \geq \mu \sqrt{k}$, 
\begin{eqnarray} \label{hz}
h(z) & \sim & \log[1+c N^{(\zeta-1)/2}(\lambda e^{(z+\mu \sqrt{k})^2/4}-1)] \\ \nonumber
& & \quad - \log[1+c N^{(\zeta-1)/2}(\lambda e^{z^2/4}-1)] \\ \nonumber
& \geq & [c \lambda+o(1)] N^{(\zeta-1)/2} e^{\mu^2 k}.
\end{eqnarray}
Since $P \{ Y_n \geq \mu \sqrt{k} \} \geq C e^{-\mu^2 k/2}/\sqrt{\log N}$,
by (\ref{hz}) and LLN,
\begin{eqnarray*}
H_N & \geq & [C+o_p(1)] N^{1-\beta+(\zeta-1)/2} e^{\mu^2 k/2}/\sqrt{\log N} \cr
& \geq & [C+o_p(1)] N^{1-\beta+(\zeta-1)/2+\rho(\beta,\zeta)+\epsilon}/\sqrt{\log N}
\end{eqnarray*}
for some $\epsilon > 0$,
and (\ref{7.8}) holds because $1-\beta+\frac{\zeta-1}{2}+\rho(\beta,\zeta)=\zeta$.

\smallskip (b) $\frac{3(1-\zeta)}{4} < \beta < 1-\zeta$ and $\rho(\beta,\zeta)=(x-y)^2$, 
where $x=\sqrt{1-\zeta}$, $y=\sqrt{1-\zeta-\beta}$.
Since $\mu^2 k=2(1+\delta)(x-y)^2 \log N+O(1)$,
by the first relation in (\ref{hz}),
$h(z) \geq C$ for $z \geq \sqrt{2(y^2-\epsilon) \log N}$ with $\epsilon > 0$ small.
Since $P \{ Y_n \geq \sqrt{2(y^2-\epsilon) \log N} \} \geq C N^{-y^2+\epsilon}/\sqrt{\log N}$, 
by LLN,
$$H_N \geq [C+o_p(1)] N^{1-\beta-y^2+ \epsilon}/\sqrt{\log N},
$$
and (\ref{7.8}) holds because $1-\beta-y^2=\zeta$. 

\medskip Case 2: $\beta>1-\zeta$ and $k=M N^{\beta+\zeta-1}$.
By the first relation in (\ref{hz}),
for $z \geq 0$,
$$h(z) \geq \tfrac{k \mu^2}{4}+O(\log N) = [\tfrac{M \mu^2}{4}+o(1)] N^{\beta+\zeta-1}.
$$
By LLN, 
$H_N \geq [\frac{M \mu^2}{8}+o_p(1)] N^\zeta$, 
and (\ref{7.8}) holds for $M > 32 \mu^{-2}$. 

\section{Proof of Theorem \ref{thm3}}

In Lemma \ref{lem5} below we provide an upper bound of the detection threshold  of the extended Mei's stopping rule,
and follow this with conditions under which this bound is exceeded under $P_\nu$.
We complete the proof by checking these conditions for various cases.
Let 
\begin{eqnarray*}
g_0(x) = g_M(x)-u \mbox{ where } g_M(x) & = & \log[1+p_0(\lambda_M e^{x/2}-1)], \cr
u & = & \log[1+p_0(\lambda_M \xi-1)],
\end{eqnarray*}
and $\xi=\lim_{t \rightarrow \infty} E e^{R_{nt}/2}$.

\begin{lem} \label{lem5}
Consider stopping rule $T_{\rm Mei}(p_0)$,
$0 < p_0 \leq 1$.
If threshold $b=Nu+\log(4 \gamma)$, 
then $E_\infty T_{\rm Mei}(p_0) \geq \gamma$.
\end{lem}

{\sc Proof}. 
If $b=Nu+\log(4 \gamma)$,
then 
$$T_{\rm Mei}(p_0) = \inf \Big\{ t: \sum_{n=1}^N g_0(R_{nt}) \geq \log(4 \gamma) \Big\}.
$$
Let $S_j = \sum_{i=1}^j Y_i$ with $Y_i$ i.i.d. N(0,1), 
and let
$$R = \sup_{j \geq 0} (\mu_0 S_j-j \mu_0^2/2). 
$$
Let $R_1, \ldots, R_N$ be an i.i.d. sample with the distribution of $R$.
Let $1 \leq t < 2 \gamma$.
Since $R_{nt}$ is bounded stochastically by $R_n$,
it follows from $E e^{g_0(R_n)}=1$,
a change of measure argument and $g_0$ monotone that
$$P_\infty \Big\{ \sum_{n=1}^N g_0(R_{nt}) \geq \log (4 \gamma) \Big\} \leq P \Big\{ \sum_{n=1}^N g_0(R_n) \geq \log (4 \gamma) \Big\} \leq (4 \gamma)^{-1}.
$$
Therefore $\{ T_{\rm Mei}(p_0) < 2 \gamma \}$ is a union of no more than $2 \gamma$ events, 
each with probability bounded by $(4 \gamma)^{-1}$ under $P_\infty$.
We conclude that $P_\infty \{ T_{\rm Mei}(p_0) < 2 \gamma \} \leq \frac{1}{2}$.
Hence $E_\infty T_{\rm Mei}(p_0) \geq \gamma$. $\wbox$

\medskip
Let $\nu \geq K_n$ and $t=\nu+k-1$,
where $k = \lfloor (1+\delta) 2 \mu^{-2} \rho(\beta,\zeta) \log N \rfloor$ for $\delta > 0$ small.
Let $U_n = \mu_0 \sum_{i=\nu}^t X_{ni}-k \mu_0^2/2$ ($\leq R_{nt}$).
Under $P_\nu$,
$U_n \sim \mbox{N}(k \mu_0 (\mu-\frac{\mu_0}{2}),k \mu_0^2)$ when $n \in {\cal N}$.
Theorem \ref{thm3} follows from 
\begin{eqnarray} \label{PgR}
\inf_{m \in J_N} P_\nu \Big\{ \sum_{n \not\in \cN} g_0(R_{nt}) \geq -N^{\zeta+\epsilon} \Big| \# \cN = m \Big\} & \rightarrow & 1, \\ \label{PgY}
\inf_{m \in J_N} P_\nu \Big\{ \sum_{n \in \cN} g_0(U_n^+) \geq 2N^{\zeta+\epsilon} \Big| \# \cN = m \Big\} & \rightarrow & 1, 
\end{eqnarray}
for some $\epsilon > 0$,
with $m \sim N^{1-\beta}$ uniformly over $m \in J_N$,
see (\ref{JN}).

The following lemma provides the framework for showing (\ref{PgR}) and (\ref{PgY}).
Let $\wtd g_0(y)=g_0(y) {\bf I}_{\{ y \leq v_N \}}$,
where 
\begin{equation} \label{vN}
v_N = (1-\zeta) \log N-\log \log N.
\end{equation}

\begin{lem} \label{lem6}
If $t \geq 4 \mu_0^{-2} (1-\zeta) \log N$ and $n \not\in \cN$,
then for all $\epsilon > 0$,
\begin{equation} \label{ER2}
E_\infty (e^{R_{nt}/2} {\bf I}_{\{ R_{nt} \leq v_N \}})  = \xi + o(N^{(\zeta-1)/2+\epsilon}).
\end{equation}
Moreover if $p_0 \sim cN^{(\zeta-1)/2}$ with $c>0$,
then
\begin{eqnarray} \label{E2}
[\inf_{y \geq 0} \wtd g_0(y)]^2 & = & O(N^{\zeta-1}), \\ \label{E3}
\sup_{y \geq 0} \wtd g_0(y) & = & O(1).
\end{eqnarray}
\end{lem}

{\sc Proof}. 
The relation (\ref{E2}) follows from 
$$| \inf_{y \geq 0} \wtd g_0(y)|=| g_0(0)| = O(p_0) = O(N^{(\zeta-1)/2}),
$$
whereas (\ref{E3}) follows from $\sup_{y \geq 0} \wtd g_0(y) = g_0(v_N) = O(1)$. 

By (\ref{Rnt}),
we can express 
$$R_{nt} = \sup_{1 \leq s \leq t} [\mu_0 S_{nst} -(t-s+1) \mu_0^2/2]^+.
$$
Extend $\{ X_{nu}: u \geq 1 \}$ to $\{ X_{nu}: -\infty<u<\infty \}$ by letting $X_{nu}$ i.i.d. N(0,1) under $P_\infty$ for $u \leq 0$.
Fix $t$ and let
$$R_n^* = \sup_{-\infty < s \leq t} [\mu_0 S_{nst}-(t-s+1) \mu_0^2/2]^+,
$$
extending the definition of $S_{nst} = \sum_{i=s}^t X_{ni}$ to $s \leq 0$.

Since $\xi = \lim_{t \rightarrow \infty} E_\infty e^{R_{nt}/2}$,
to show (\ref{ER2}), 
it suffices to show that 
\begin{eqnarray} \label{R1}
E_\infty (e^{R_{nt}/2} {\bf I}_{\{ R_{nt} > v_N \}}) & = & o(N^{(\zeta-1)/2+\epsilon}), \\ \label{R2}
E_\infty (e^{R^*_n/2} {\bf I}_{\{ R^*_n > R_{nt} \}}) & = & o(N^{(\zeta-1)/2+\epsilon}).
\end{eqnarray}

We conclude (\ref{R1}) from (\ref{limR}) and (\ref{vN}). 
Let $Q = \sup_{j \geq t} (\mu_0 S_j-j \mu_0^2/2)$ and $R' = \sup_{j \geq 0} (\omega S_j-j \omega^2 \mu_0^2/2)$ for some $\omega > \frac{1}{2}$. 
By (\ref{limR}), 
for $x \geq 0$,
\begin{eqnarray*}
P_\infty \{ R^*_n > R_{nt}, R^*_n \geq x \} \leq P \{ Q \geq x \} & \leq & P \{ R' \geq \omega x+t (\omega-\omega^2) \mu_0^2/2 \} \cr
& = & O(e^{-\omega x-t (\omega-\omega^2)\mu_0^2/2}).
\end{eqnarray*}
Hence by selecting $\omega$ close enough to $\frac{1}{2}$,
it follows that
$$E_\infty (e^{R^*_n/2} {\bf I}_{\{ R^*_n > R_{nt} \}}) = \int_{-\infty}^{\infty} \tfrac{1}{2} e^{x/2} P_\infty \{ R_n^* > R_{nt}, R_n^* \geq x \} dx = O(e^{-t \mu_0^2/8+\epsilon}), 
$$
and (\ref{R2}) holds for $t \geq 4 \mu_0^{-2} (1-\zeta) \log N$. $\wbox$  

\medskip
We conclude (\ref{PgR}) from (\ref{ER2}), 
$p_0 \sim cN^{(\zeta-1)/2}$ and LLN. 
We note that indeed $t (=\nu+k-1) \geq 4 \mu_0^{-2} (1-\zeta) \log N$ when

(a) $\frac{1-\zeta}{2} < \beta \leq \frac{3(1-\zeta)}{4}$, $\rho(\beta,\zeta)=\beta-\frac{1-\zeta}{2}$, $\mu_0=2 \mu$, $\nu \geq 2 \mu^{-2} (1-\zeta-\beta) \log N$,

(b) $\frac{3(1-\zeta)}{4} < \beta < 1-\zeta$, $\mu_0 = \mu \sqrt{\frac{1-\zeta}{\rho(\beta,\zeta)}}$, $\nu \geq 2 \mu^{-2} \rho(\beta,\zeta) \log N$.

\noindent It remains for us to check (\ref{PgY}) on: 

\medskip
(a) $\frac{1-\zeta}{2} < \beta < \frac{3(1-\zeta)}{4}$.
Since $\mu_0=2 \mu$,
we have $U_n \sim$ N(0,$k \mu_0^2$) when $n \in \cN$.
Hence 
\begin{eqnarray} \label{Eg1}
& & E_\nu [\wtd g_0(U_n^+)|n \in \cN] \\ \nonumber 
& & \qquad \sim p_0 [E_\nu (e^{U_n^+/2} {\bf I}_{\{ U_n \leq v_N \}}|n \in \cN)-\xi] \\ \nonumber
& & \qquad \sim cN^{(\zeta-1)/2} e^{k \mu_0^2/8} \int_{-\infty}^{v_N} \tfrac{1}{\sqrt{2 k \pi \mu_0^2}} e^{-(y-k \mu_0^2/2)^2/(2k \mu_0^2)} dy \\ \nonumber
& & \qquad = cN^{(\zeta-1)/2} e^{k \mu_0^2/8} \Phi(\tfrac{v_N-k \mu_0^2/2}{\mu_0 \sqrt{k}}).
\end{eqnarray}
Check that $e^{k \mu_0^2/8} = N^{[1+\delta+o(1)] \rho(\beta,\zeta)} \geq N^{\beta+(\zeta-1)/2+2 \epsilon}$ for $\epsilon > 0$ small and $N$ large.
Moreover $\rho(\beta,\zeta) < \frac{1-\zeta}{4}$,
therefore $v_N > \frac{k \mu_0^2}{2} [\sim 4(1+\delta) \rho(\beta,\zeta) \log N]$ for $\delta > 0$ small.
Hence by (\ref{Eg1}),
\begin{equation} \label{Eplus}
E_\nu [\wtd g_0(U_n^+)|n \in \cN] \sim [c+o(1)] N^{\beta+\zeta-1+2 \epsilon}.
\end{equation}
By (\ref{E2}), (\ref{E3}) and (\ref{Eplus}), 
we can conclude 
$$E_\nu [\wtd g_0^2(U_n^+)|n \in \cN] = O(|E_\nu [\wtd g_0(U_n^+)|n \in \cN]|),
$$
and (\ref{PgY}) then follows from (\ref{Eplus}), 
Chebyshev's inequality and $g_0 \geq \wtd g_0$.

\medskip
(b) $\frac{3(1-\zeta)}{4} \leq \beta < 1-\zeta$.
For $n \in \cN$,
express $U_n = k \mu_0 (\mu-\frac{\mu_0}{2})+\sqrt{k} \mu_0 Y_n$,
with $Y_n \sim$ N(0,1).
Let $\cN_1 = \{ n \in \cN: Y_n \geq \sqrt{2(1-\zeta-\beta-2 \epsilon) \log N} \}$ for $\epsilon > 0$ satisfying 
\begin{equation} \label{zb}
1-\zeta-\beta-2 \epsilon \geq (1-\zeta-\beta)/(1+\delta).
\end{equation}
By LLN,
\begin{eqnarray} \label{hN1}
\# \cN_1 & = & (\# \cN) [C+o_p(1)] N^{\zeta+\beta-1+2 \epsilon}/\sqrt{\log N} \\ \nonumber
& = & [C+o_p(1)] N^{\zeta+2 \epsilon}/\sqrt{\log N}.
\end{eqnarray}
Let $r=\mu_0/\mu (=\sqrt{\frac{1-\zeta}{\rho(\beta,\zeta)}} \leq 2)$.
Since $r-1=\sqrt{\frac{1-\zeta-\beta}{\rho(\beta,\zeta)}}$,
by (\ref{k}) and (\ref{zb}),
for $n \in \cN_1$ with $N$ large,
\begin{eqnarray*}
U_n & \geq & k \mu^2(r-\tfrac{r^2}{2})+\mu r(r-1) \sqrt{\tfrac{2k \rho(\beta,\zeta) \log N}{1+\delta}} \cr
& = & 2 \rho(\beta,\zeta) \log N[(1+\delta)(r-\tfrac{r^2}{2})+r^2-r]+O(1) \cr
& \geq & r^2 \rho(\beta,\zeta) \log N = (1-\zeta) \log N, \cr
g_0(U_n^+) & \geq & \log[1+p_0(N^{(\zeta-1)/2}-1)] - \log[1+p_0(\xi-1)] \rightarrow \log(1+c).
\end{eqnarray*}
Hence by (\ref{hN1}),
$$\sum_{n \in \cN_1} g_0(U_n^+) \geq [C+o_p(1)] N^{\zeta+2 \epsilon}/\sqrt{\log N}.
$$
This,
combined with
\begin{eqnarray*}
\sum_{n \in \cN \setminus \cN_1} g_0(U_n^+) & \geq & -[C+o_p(1)] N^{1-\beta} \log p_0 \cr
& \sim &- [C+o_p(1)] N^{1-\beta+(\zeta-1)/2},
\end{eqnarray*}
and noting that $1-\beta+\frac{\zeta-1}{2} \leq 1-\frac{5}{4}(1-\zeta) < \zeta$,
shows (\ref{PgY}). 

\medskip
{\bf Acknowledgments}. 
We thank an Associate Editor and the referees for their insights and helpful comments.
The appendices below are due to questions and suggestions from the referees.

\begin{appendix}
\section{Motivations behind (\ref{k=})}

In view of the need to satisfy (\ref{p0}),
we choose $k$ to be the ``largest" possible such that 
\begin{equation} \label{A1}
E_1 \ell_{\bullet 1k} < 2N^\zeta/3.
\end{equation}
Under $P_1$,
$Z_{nk} \sim (1-p)$N(0,1)$+p$N($\mu \sqrt{k},1$).
Let $Y \sim$ N(0,1).
Since
\begin{equation} \label{A3}
\log(1+x) \leq x, \qquad E e^{Y \mu \sqrt{k}}=e^{k \mu^2/2},
\end{equation}
it follows that
\begin{eqnarray} \label{A2}
E_1 \ell_{\bullet 1k} & \{ = & N(1-p) E \log[1+p(e^{Y \mu \sqrt{k}-k \mu^2/2}-1)] \\ \nonumber
& & \quad +Np E \log [1+p(e^{Y \mu \sqrt{k}+k \mu^2/2}-1)] \} \\ \nonumber
& \leq & Np E \log[1+p(e^{Y \mu \sqrt{k}+k \mu^2/2}-1)].
\end{eqnarray}

Case 1(a): $\frac{1-\zeta}{2} < \beta \leq \frac{3(1-\zeta)}{4}$.
It follows from applying (\ref{A3}) on (\ref{A2}) that
\begin{equation} \label{A4}
E_1 \ell_{\bullet 1k} \leq Np^2 (e^{k \mu^2}-1) \sim N^{1-2 \beta+(k \mu^2/\log N)}.
\end{equation}
Hence choosing $k = \lfloor (1-\delta) \mu^{-2} (2 \beta+\zeta-1) \log N \rfloor$ as in (\ref{k=}) ensures $E_1 \ell_{\bullet 1k} = o(N^\zeta)$,
and so (\ref{A1}) holds.

\medskip
Case 1(b): $\frac{3(1-\zeta)}{4} < \beta < 1-\zeta$.
The inequality in (\ref{A4}) is further sharpened to allow for larger $k$ satisfying (\ref{A1}).
Let $\omega$ be the root of 
\begin{equation} \label{A5}
e^{\omega \mu \sqrt{k}+k \mu^2/2} = N^\beta (\sim p^{-1}).
\end{equation}
By (A5) applying the inequalities 
\begin{equation} \label{Alog}
\log(1+x) \leq \left\{ \begin{array}{ll} x & \mbox{ if } -1 < x < 1, \cr
\log 2+\log x & \mbox{ if } x \geq 1, 
\end{array} \right.
\end{equation}
on (\ref{A2}) results in 
\begin{eqnarray} \label{A6}
E_1 \ell_{\bullet 1k} & \leq & Np^2 \int_{-\infty}^\omega \tfrac{1}{\sqrt{2 \pi}} e^{-z^2/2+z \mu \sqrt{k}+k \mu^2/2} dz + O(Npk e^{-\omega^2/2}) \\ \nonumber
& = & Np^2 e^{k \mu^2} \Phi(\omega-\mu \sqrt{k})+O(Npk e^{-\omega^2/2}) \\ \nonumber
& = & O(Npk e^{-\omega^2/2}).
\end{eqnarray}
By (\ref{A5}),
\begin{equation} \label{A7}
\omega \mu \sqrt{k}+k \mu^2/2 = \beta \log N (\Rightarrow \mu \sqrt{k} = -\omega+\sqrt{\omega^2+2 \beta \log N}), 
\end{equation}
and by (\ref{A6}),
we satisfy (\ref{A1}) if 
\begin{equation} \label{A8}
(1-\beta) \log N-\omega^2/2 < \zeta \log N (\Rightarrow \omega > \sqrt{2(1-\beta-\zeta) \log N}).
\end{equation}
Combining (\ref{A7}) and (\ref{A8}) leads to $k < 2 \mu^{-2} (\sqrt{1-\zeta}-\sqrt{1-\zeta-\beta})^2 \log N$. 
Hence the choice of $k = \lfloor (1-\delta) 2 \mu^{-2} (\sqrt{1-\zeta}-\sqrt{(1-\zeta-\beta})^2 \log N \rfloor$ in (\ref{k=}).

\smallskip
Case 2: $\beta > 1-\zeta$.
By (\ref{A2}) and (\ref{Alog}),
choosing $k= \lfloor \delta N^{\beta+\zeta-1} \rfloor$ as in (\ref{k=}) ensures that 
$$E_1 \ell_{\bullet 1k}\leq [1+o(1)] Np E(Y \mu \sqrt{k}+k \mu^2/2) \sim \delta \mu^2 N^\zeta/2,
$$
and (\ref{A1}) indeed holds for $\delta > 0$ small.

\section{Minimum detection delay under the minimax setting}

Let ${\bf I}_N = (I_1, \ldots, I_N)$, 
where $I_n = {\bf I}_{\{ n \in \cN \}}$.
Let $E_{\nu,{\bf I}_N}$ denote expectation with respect to $X_{nt} \sim$ N($\mu_{nt},1$),
with $\mu_{nt} = \mu I_n {\bf I}_{\{ t \geq \nu \}}$.
For a given stopping rule $T$,
define 
$$D_{N,m}(T) = \sup_{1 \leq \nu < \infty} \Big[ \max_{{\bf I}_N: \sum I_n = m} E_{\nu,{\bf I}_N} (T-\nu+1|T \geq \nu) \Big].
$$ 
The following is an analogue of Theorem \ref{thm1} on a minimax setting.

\begin{thm} \label{thmB}
Let $T$ be a stopping rule such that {\rm ARL}$(T) \geq \gamma$,
with $\log \gamma \sim N^\zeta$ for some $\zeta > 0$.
Let $m \sim N^{1-\beta}$ for some $0 < \beta < 1$.

\smallskip
{\rm (a)} If $\tfrac{1-\zeta}{2} < \beta < 1-\zeta$,
then
$$\liminf_{N \rightarrow \infty} \frac{D_{N,m}(T)}{\log N} \geq 2 \mu^{-2} \rho(\beta,\zeta).
$$

{\rm (b)} If $\beta > 1-\zeta$, 
then
$$\liminf_{N \rightarrow \infty} \frac{\log D_{N,m}(T)}{\log N} \geq \beta+\zeta-1.
$$
\end{thm}

\medskip
{\sc Proof}. 
Let $k$ be chosen as in (\ref{k=}).
By Lemma \ref{lem1} we can find $s \geq 1$ such that
\begin{equation} \label{B1}
P_\infty \{ T \geq s+k | T \geq s \} \geq 1- k/\gamma.
\end{equation} 
Let $t=s+k-1$, 
and consider the test,
conditional on $T \geq s$,
of
$$\begin{array}{rl}
H_0: & X_{nu} \sim \mbox{ N(0,1) for } 1 \leq n \leq N, 1 \leq u \leq t, \cr
\mbox{vs } H_{s,m}: & X_{nu} \sim \mbox{N}(\mu {\bf I}_{\{ u \geq s, n \in \cN \}},1) \mbox{ for  } 1 \leq n \leq N, 1 \leq u \leq t, \cr
& \mbox{ with } \cN \mbox{ a random subset of } \{ 1, \ldots, N \} \mbox{ of size } m. 
\end{array}
$$
By (\ref{B1}) the test rejecting $H_0$ when $T <s+k$ has Type I error probability not exceeding $k/\gamma$.

Let $\cA_j = \{ \cN: \# \cN=j \}$.
At time $t$, 
the (conditional) likelihood ratio between $H_{s,m}$ and $H_0$ is $L_m (=L_{mst})$, 
where
$$L_j = {N \choose j}^{-1} \sum_{\cN \in \cA_j} \Big( \prod_{n \in \cN} e^{Z_n \mu \sqrt{k}-k \mu^2/2} \Big),
\quad Z_n = Z_{nst}.
$$
Let $P_{s,m}$ ($E_{s,m}$) denote probability (expectation) with respect to $H_{s,m}$.

We shall check on various cases below that
\begin{equation} \label{Pss}
P_{s,m} \{ L_m \geq J \} \rightarrow 0, \quad J=\exp(2N^\zeta/3).
\end{equation}
Let $B$ be such that $P_{s,m} \{ J \geq L_m \geq B \} = \exp(-N^\zeta/4)$.
It follows from (\ref{Pss}) that
\begin{equation} \label{B3}
P_{s,m} \{ L_m \geq B \} (= P_{s,m} \{ L_m \geq B| T \geq s \}) \rightarrow 0, 
\end{equation}
and that for $N$ large,
\begin{eqnarray}
\label{B4}
& & P_{\infty} \{ L_m \geq B \} (= P_\infty \{ L_m \geq B | T \geq s \}) \geq P_\infty \{ J \geq L_m \geq B \} \\ \nonumber
& = & E_{s,m} (L_m^{-1} {\bf I}_{\{ J \geq L_m \geq B \}}) \geq J^{-1} \exp(-N^\zeta/4) \geq k/\gamma.
\end{eqnarray}
By (\ref{B1}), (\ref{B4}) and the Neyman-Pearson Lemma,
the test rejecting $H_0$ when $L_m \geq B$ is at least as powerful as the one based on $T$,
that is
\begin{equation} \label{B5}
P_{s,m} \{ T \geq s+k |T \geq s \} \geq P_{s,m} \{ L_m < B \}.
\end{equation}
It follows from (\ref{B3}) and (\ref{B5}) that
$$D_{N,m}(T) \geq E_{s,m} (T-s+1|T \geq s) \geq k P_{s,m} \{ T \geq s+k|T \geq s \} = k[1+o(1)],
$$
and the proof of Theorem \ref{thmB} is complete. 
$\wbox$

\medskip
We shall now proceed to check (\ref{Pss}).
Let $p_1 = 2N^{-\beta}$ and 
\begin{equation} \label{Lp1}
L(p_1) = \prod_{n=1}^N (1-p_1+p_1 e^{Z_n \mu \sqrt{k}-k \mu^2/2}) \Big[ = \sum_{j=0}^N (1-p_1)^{N-j} p_1^j {N \choose j} L_j \Big].
\end{equation}
Since $Z_n \sim$ N($\mu \sqrt{k},1$) if $n \in \cN$ and $Z_n \sim$ N(0,1) if $n \not\in \cN$,
it follows that
$$E e^{Z_n \mu \sqrt{k}-k \mu^2/2}= \left\{ \begin{array}{ll} e^{k \mu^2} & \mbox{ if } n \in \cN, \cr
1 & \mbox{ if } n \not\in \cN. \end{array} \right.
$$
Therefore by (\ref{Lp1}),
\begin{equation} \label{Ess1}
E_{s,m} L(p_1) = (1-p_1+p_1 e^{k \mu^2})^m,
\end{equation}
the exponent $m$ in (\ref{Ess1}) due to $\# \cN=m$ for each $\cN$ under $H_{s,m}$.
By the monotonicity $E_{s,m} L_1 \leq \cdots \leq E_{s,m} L_N$, 
and by $P \{ W \geq m \} \rightarrow 1$ for $W \sim$ Binomial($N,p_1$),
it follows from (\ref{Lp1}) that
\begin{equation} \label{Ess2}
E_{s,m} L(p_1) \geq P \{ W \geq m \} E_{s,m} L_m = [1+o(1)] E_{s,m} L_m.
\end{equation}
By (\ref{Ess1}), (\ref{Ess2}) and Markov's inequality, 
to show (\ref{Pss}) it suffices to show that
\begin{equation} \label{show}
(1-p_1+p_1 e^{k \mu^2})^m = o(\exp(2N^\zeta/3)),
\end{equation} 
and this can be easily done for the following cases.

\medskip
Case 1(a): $\tfrac{1-\zeta}{2} < \beta \leq \tfrac{3(1-\zeta)}{4}$, 
$k = \lfloor (1-\delta) \mu^{-2} (2 \beta+\zeta-1) \log N \rfloor$.
We show (\ref{show}) by applying the inequality
$$(1-p_1 + p_1 e^{k \mu^2})^m \leq \exp(mp_1 e^{k \mu^2}).
$$

\smallskip
Case 2: $\beta > 1-\zeta$, $k = \lfloor \delta N^{\beta+\zeta-1} \rfloor$, $\delta > 0$ small.
We show (\ref{show}) by applying the inequalities (for large $N$),
$$(1-p_1 +p_1 e^{k \mu^2})^m \leq (2p_1 e^{k \mu^2})^m \leq e^{k \mu^2 m}.
$$

\smallskip
The final case below is more complicated.
Additional truncation arguments are needed to show (\ref{Pss}).

\smallskip
Case 1(b): $\tfrac{3(1-\zeta)}{4} < \beta < 1-\zeta$,
$k = \lfloor (1-\delta) 2 \mu^{-2} (x-y)^2 \log N \rfloor$,
where $x= \sqrt{1-\zeta}$ and $y=\sqrt{1-\zeta-\beta}$.
The outline of the arguments needed to show (\ref{Pss}) is as follows.

\begin{enumerate}
\item Let $\wtd Z_n = \min(Z_n, \omega)$,
where
$$\omega (=\omega_N) = \sqrt{2(1-\zeta) \log N+2 \log \log N} (\doteq x \sqrt{2 \log N}).
$$
Let $p_1 = 2N^{-\beta}$ and 
\begin{equation} \label{Ltp1}
\wtd L(p_1) = \prod_{n=1}^N (1-p_1+p_1 e^{\tilde Z_n \mu \sqrt{k}-k \mu^2/2}).
\end{equation}
Show that $E_{s,m} \wtd L(p_1) = o(J^{1/2}) [=o(\exp(N^{\zeta}/3))]$.

\item Argue that we have monotonicity $E_{s,m} \wtd L_1 \leq \cdots \leq E_{s,m} \wtd L_N$,
where
$$\wtd L_j = {N \choose j}^{-1} \sum_{\cN \in \cA_j} \Big( \prod_{n \in \cN} e^{\tilde Z_n \mu \sqrt{k}-k \mu^2/2} \Big),
$$
and conclude that
\begin{equation} \label{Esm}
E_{s,m} \wtd L(p_1) \geq P \{ W \geq m \} E_{s,m} \wtd L_m = [1+o(1)] E_{s,m} \wtd L_m,
\end{equation}
where $W \sim$ Binomial($N,p_1$).

\item Let $C>0$ and $\wht L_m = L_m {\bf I}_G$,
where $G(=G_N$) is the event that
$$\max_{1 \leq n \leq N} Z_n \leq C \sqrt{\log N}, \quad F_N := \# \{ n: Z_n > \omega \} \leq N^\zeta/(\log N)^{5/4}.
$$
Show that uniformly under $G$,
$$\max_{\cN \in \cA_m} \Big( \prod_{n \in \cN} e^{(Z_n - \tilde Z_n) \mu \sqrt{k}} \Big) = o(J^{1/2}) [= o(\exp(N^\zeta/3))],
$$
and conclude that $\wht L_m/\wtd L_m = o(J^{1/2})$.

\item Show that for $C$ large,
$P_{s,m}(G_N) \rightarrow 1$ and so $P_{s,m} \{ L_m > \wht L_m \} \rightarrow 0$.
\end{enumerate}
By steps 1, 2 and Markov's inequality, 
$P_{s,m} \{ \wtd L_m \geq J^{1/2} \} \rightarrow 0$.
By step~3 we can further conclude that $P_{s,m} \{ \wht L_m \geq J \} \rightarrow 0$,
and (\ref{Pss}) then follows from step~4.
We shall now provide details to the above outline.

\begin{enumerate}
\item If $n \not\in \cN$, then $E e^{\tilde Z_n \mu \sqrt{k}-k \mu^2/2} \leq 1$,
and if $n \in \cN$,
then
\begin{eqnarray*}
E e^{\tilde Z_n \mu \sqrt{k}-k \mu^2/2} & = & e^{k \mu^2} \Phi(\omega-2 \mu \sqrt{k})+[1-\Phi(\omega - \mu \sqrt{k})] e^{\omega \mu \sqrt{k}-k \mu^2/2} \cr
& = & o(N^{2(x-y)^2-(2y-x)^2})+o(N^{-y^2+2x(x-y)-(x-y)^2}) \cr
& = & o(N^{x^2-2y^2}).
\end{eqnarray*}
Since $\# \cN=m$ for each $\cN$ under $H_{s,m}$,
by (\ref{Ltp1}),
$$E_{s,m} \wtd L(p_1) \leq [1+p_1 o(N^{x^2-2y^2})]^m \leq \exp[mp_1 o(N^{x^2-2y^2})]=o(J^{1/2}).
$$

\item The monotonicity follows from $\wtd Z_n$ stochastically larger when $n \in \cN$ compared to when $n \not\in \cN$,
whereas the inequality in (\ref{Esm}) follows from the monotonicity and the expansion
$$\wtd L(p_1) = \sum_{j=0}^N (1-p_1)^{N-j} p_1^j {N \choose j} \wtd L_j.
$$

\item Under $G$, there exists $\wtd C>0$ not depending on $N$ such that for all $\cN \in \cA_m$,
\begin{eqnarray*}
\prod_{n \in \cN} e^{(Z_n-\tilde Z_n) \mu \sqrt{k}} & \leq & \exp(F_N C \sqrt{\log N} \mu \sqrt{k}) \cr
& \leq & \exp \Big( \tfrac{N^\zeta}{(\log N)^{5/4}} \cdot \wtd C \log N \Big) = o(J^{1/2}).
\end{eqnarray*}

\item Let $\bar \Phi(\cdot)=1-\Phi(\cdot)$. 
We apply Markov's inequality to show $P_{s,m}(G_N) \rightarrow~1$ by checking that
\begin{equation} \label{B42}
m \bar \Phi(\omega-\mu \sqrt{k})+(N-m) \bar \Phi(\omega) = o \Big( \tfrac{N^\zeta}{(\log N)^{5/4}} \Big),
\end{equation}
and that for $C$ large,
\begin{equation} \label{B41}
m \bar \Phi(C \sqrt{\log N}-\mu \sqrt{k})+(N-m) \bar \Phi(C \sqrt{\log N}) \rightarrow 0.
\end{equation}
By Mill's inequality,
(\ref{B41}) holds for $C$ large and $N \bar \Phi(\omega) = o \Big( \tfrac{N^\zeta}{(\log N)^{5/4}} \Big)$.
Moreover 
$$\limsup_{N \rightarrow \infty} \log_N [ m \bar \Phi(\omega-\mu \sqrt{k})] < 1-\beta-y^2 = \zeta,
$$
and so (\ref{B42}) holds as well.
\end{enumerate}
\end{appendix}

\end{document}